\documentclass[11pt,a4paper]{amsart}
\usepackage{graphics,array,amsmath,amssymb,epsf,graphicx}

\newtheorem{thm}{Theorem}
\newtheorem{lem}[thm]{Lemma}

\theoremstyle{definition}

\begin{document}
\title[Linearly embedded graphs]
{Linearly embedded graphs in 3-space with homotopically free exteriors}
\author{Youngsik Huh}
\author{Jung Hoon Lee}
\address{Department of Mathematics,
College of Natural Sciences, Hanyang University, Seoul 133-791,
Korea} \email{yshuh@hanyang.ac.kr}
\address{Department of Mathematics and Institute of Pure and Applied Mathematics,
Chonbuk National University,
Jeonju 561-756,
Korea} \email{junghoon@jbnu.ac.kr}




\begin{abstract}
An embedding of a graph into $\mathbb{R}^3$ is said to be {\em linear}, if any edge of the graph  is sent to be a line segment. And we say that an embedding $f$ of a graph $G$ into $\mathbb{R}^3$ is {\em free}, if $\pi_1(\mathbb{R}^3-f(G))$ is a free group.
It was known that for any complete graph its linear embedding is always free \cite{Nicholson}.

In this paper we investigate the freeness of linear embeddings considering the number of vertices.
It is shown that for any simple connected graph with at most 6 vertices, if its minimal valency is at least 3, then its linear embedding is always free. On the contrary when the number of vertices is much larger than the minimal valency or connectivity, the freeness may not be an intrinsic property of such graphs. In fact we show that for any $n \geq 1$ there are infinitely many connected graphs with minimal valency $n$ which have non-free linear embeddings, and furthermore, that there are infinitely many $n$-connected graphs which have non-free linear embeddings.
\end{abstract}

\maketitle



\section{Introduction}
Let $G$ be a finite connected graph and $f:G \rightarrow \mathbb{R}^3$ be an embedding of $G$ into the Euclidean 3-space $\mathbb{R}^3$. If the fundamental group $\pi_1(\mathbb{R}^3-f(G))$ is free, then we say that the embedding $f$ is {\em free}.

The freeness of fundamental group plays a key role in detecting the unknottedness of graphs in $\mathbb{R}^3$. A graph embedded into $\mathbb{R}^3$(or its embedding into $\mathbb{R}^3$)  is said to be {\em unknotted} if it lies on an embedded surface in $\mathbb{R}^3$ which is homeomorphic to the 2-sphere.
It is known that a simple closed curve embedded in $\mathbb{R}^3$ is unknotted if and only if it is free \cite{Pa}. This result was generalized by Scharlemann and Thompson. They  proved that for any planar graph $G$ its embedding $f$ is unknotted if and only if $\pi_1(\mathbb{R}^3-f(H))$ is free for every subgraph $H$ of $G$ \cite{ST,Wu,Go}. Furthermore this criterion is valid even when determining whether a graph has a linkless embedding. It was proved, by Robertson, Seymour and Thomas, that a graph $G$ has a linkless embedding if and only if it has an embedding $f$ such that $\pi_1(\mathbb{R}^3-f(H))$ is free for every subgraph $H$ of $G$ \cite{RST}.

\begin{figure}[h]
\centering
\includegraphics[width=8cm]{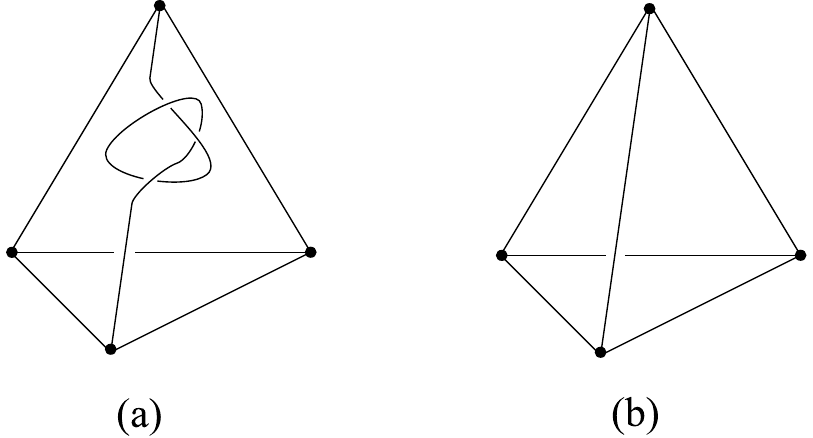}
\caption{}
\label{fig1-1}
\end{figure}

In this paper, basically we are interested in the freeness of graph in $\mathbb{R}^3$, but our viewpoint is different from the previous works in the above. We begin with observing two specific embeddings of the complete graph $K_4$. In Figure \ref{fig1-1}-(a) an edge forms a locally knotted arc, hence the embedding illustrated in the figure is not free. For any graph with a cycle, we may construct a non-free embedding in this way. On the contrary, in Figure \ref{fig1-1}-(b), every edge is a line segment. Therefore the embedded edges constitute the 1-skeleton of a tetrahedron, and the fundamental group of its complement is free in consequence. Motivated by this example, the freeness of such {\em linear} embeddings is studied in this paper.

An embedding of a graph into $\mathbb{R}^3$ will be said to be {\em linear} if each edge of the graph is sent to be a line segment. Note that, for a graph to be linearly embeddable, it should be simple, that is, have no multiple edges between any two vertices and no loop edges.

In the late 1980s, V. Nicholson proved the following theorem which says that the freeness of linear embedding is an intrinsic property of complete graphs:
\begin{thm} {\rm \cite{Nicholson}} \label{thm1}
Every linear embedding of the complete graph $K_n$ is free.
\end{thm}

\begin{figure}
\centering
\includegraphics[width=8cm]{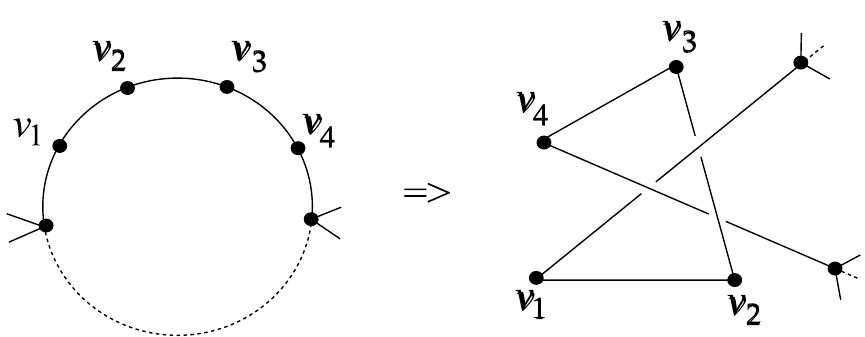}
\caption{}
\label{fig1-2}
\end{figure}

\noindent The first aim of this paper is to investigate small graphs with such property.
Let $G$ be a connected graph. Suppose that a cycle of $G$ contains $n$ consecutive vertices $v_1, \ldots, v_n$ such that the valency of $v_i$ in $G$ is two for every $i$. Then for $n \geq 4$ we can construct a linear embedding of $G$ so that the path traversing the vertices forms a locally knotted arc, and consequently the embedding is not free. See Figure \ref{fig1-2} for example. To avoid such local knottedness, the graphs in our consideration are restricted so that every vertex is of valency at least three.

\begin{figure}
\centering
\includegraphics[width=8cm]{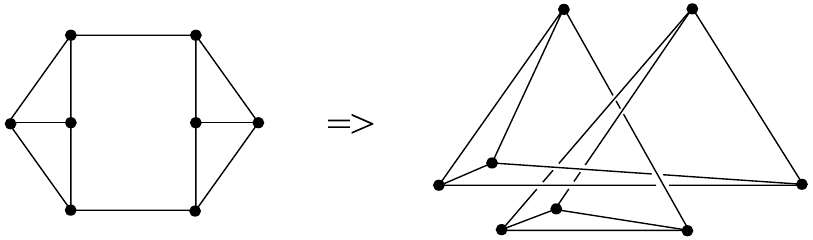}
\caption{}
\label{fig1-3}
\end{figure}

Let $V(G)$ be the set of vertices of $G$. For a vertex $v$, let $d(v)$ denote its valency. And $\delta(G)$ will denote $\mbox{min}\{d(v)|v\in V(G)\}$.  Then we have the following theorem:
\begin{thm} \label{thm2}
Let $G$ be a simple connected graph with $|V(G)| \leq 6$.
If the minimal valency $\delta(G)$ is at least $3$, then every linear embedding of $G$ is free.
\end{thm}

Here we remark that the condition, $\delta \geq 3$, is not sufficient to guarantee the freeness for graphs with more than 7 vertices. For example Figure \ref{fig1-3} shows a graph with $\delta =3$ and $|V|=8$ which has a non-free linear embedding. The embedding is obtained by adding 6 additional line segments to a hexagonal trefoil knot so that the additional edges cobound some disks with subarcs of the knot. In consequence the fundamental group of its complement is a free product of the fundamental group of the trefoil knot complement and a free group. Since the fundamental group of the trefoil knot complement is not free, our group is also not free.

Motivated by this observation we show that if the number of vertices of a graph is relatively larger than its minimal valency or connectivity, then it may have a linear embedding which is not free:
\begin{thm} \label{thm3}
For any $n \geq 1$, there are infinitely many simple connected graphs with minimal valency $n$ which have non-free linear embeddings.
\end{thm}
\begin{thm} \label{thm4}
For any $n \geq 1$, there are infinitely many $n$-connected graphs which have non-free linear embeddings.
\end{thm}
\noindent The proofs of the two theorems are constructive. In Theorem \ref{thm3} the constructed graph has at least $6(n+1)$ vertices, and in Theorem \ref{thm4}, at least $12n$ vertices. Note that for the complete graph $K_n$, the number of vertices, valency and connectivity are $n$, $n-1$ and $n$, respectively.

The rest of this paper is devoted to proving the theorems. Theorems \ref{thm3} and \ref{thm4} are proved in Sections 2 and 3, respectively. The proof of Theorem \ref{thm2} is given in the final section for our convenience.

\section{Proof of Theorem \ref{thm3}}
\begin{figure}[ht!]
\centering
\includegraphics[width=5.5cm]{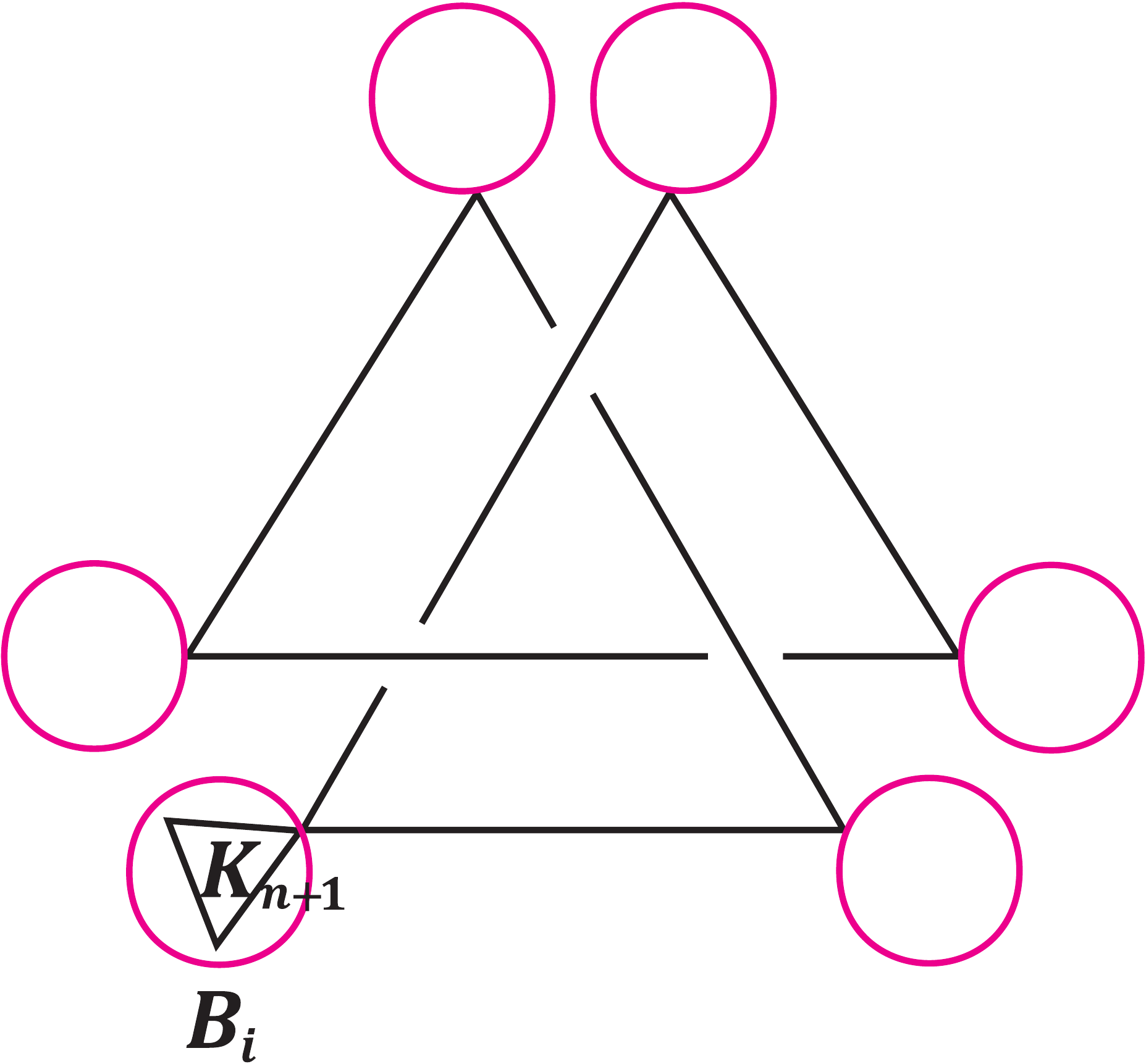}
\caption{A linear embedding of $G$ which is not free}\label{fig2-1}
\end{figure}

Let $C_6$ be the cycle graph with only six vertices $v_1, v_2, \ldots , v_6$.
And let $G_i$ be a copy of the complete graph $K_{n+1}$ for $1 \leq i \leq 6$.
Identifying each $v_i$ with a vertex of $G_i$, we obtain a simple connected graph $G$ with
$6(n+1)$ vertices and $\delta (G) = n$.
And then construct a linear embedding $f$ of $G$  into $\mathbb{R}^3$ as follows.
Figure \ref{fig2-1} depicts the embedding.

\begin{itemize}
\item Embed $C_6$ so that $f(C_6)$ is a hexagonal trefoil knot.
\item Take mutually disjoint small 3-balls $B_1, B_2, \ldots , B_6$ so that each $B_i$ meets $f(C_6)$ only at the vertex $f(v_i)$.
\item Embed each $G_i$ linearly into $B_i$ so that it meets the boundary 2-sphere $S_i=\partial B_i$ only at the identified vertex $f(v_i)$.
\end{itemize}

\begin{figure}[ht!]
\centering
\includegraphics[width=8.5cm]{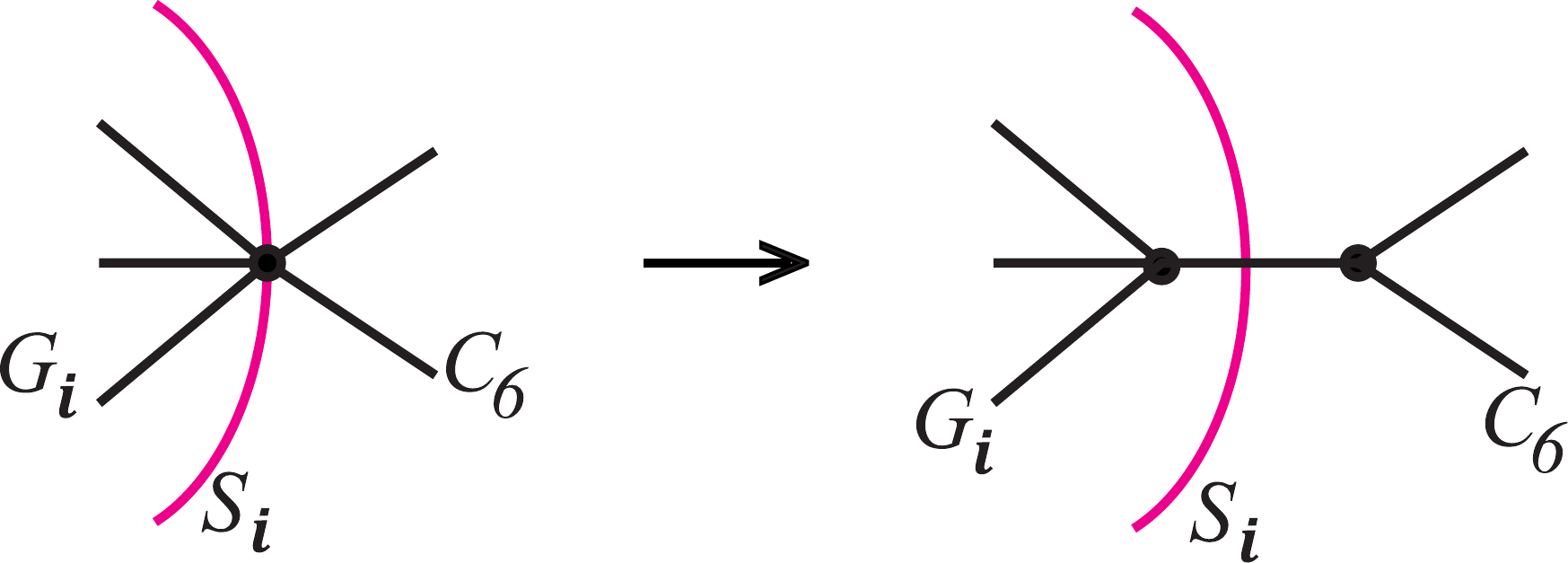}
\caption{Subdivision at $v_i$ }\label{fig2-2-0}
\end{figure}

Now we show that $\pi_1(\mathbb{R}^3-f(G))$ is not free. For our convenience let $G$ denote the embedded graph $f(G)$ itself.
Firstly modify $G$ slightly by subdividing each vertex $v_i$ as seen in Figure \ref{fig2-2-0}.
Then the complement of $G$ has the same homotopy type as that of the new graph. Again for our convenience let $G$ denote the new graph.

\begin{figure}[ht!]
\centering
\includegraphics[width=6cm]{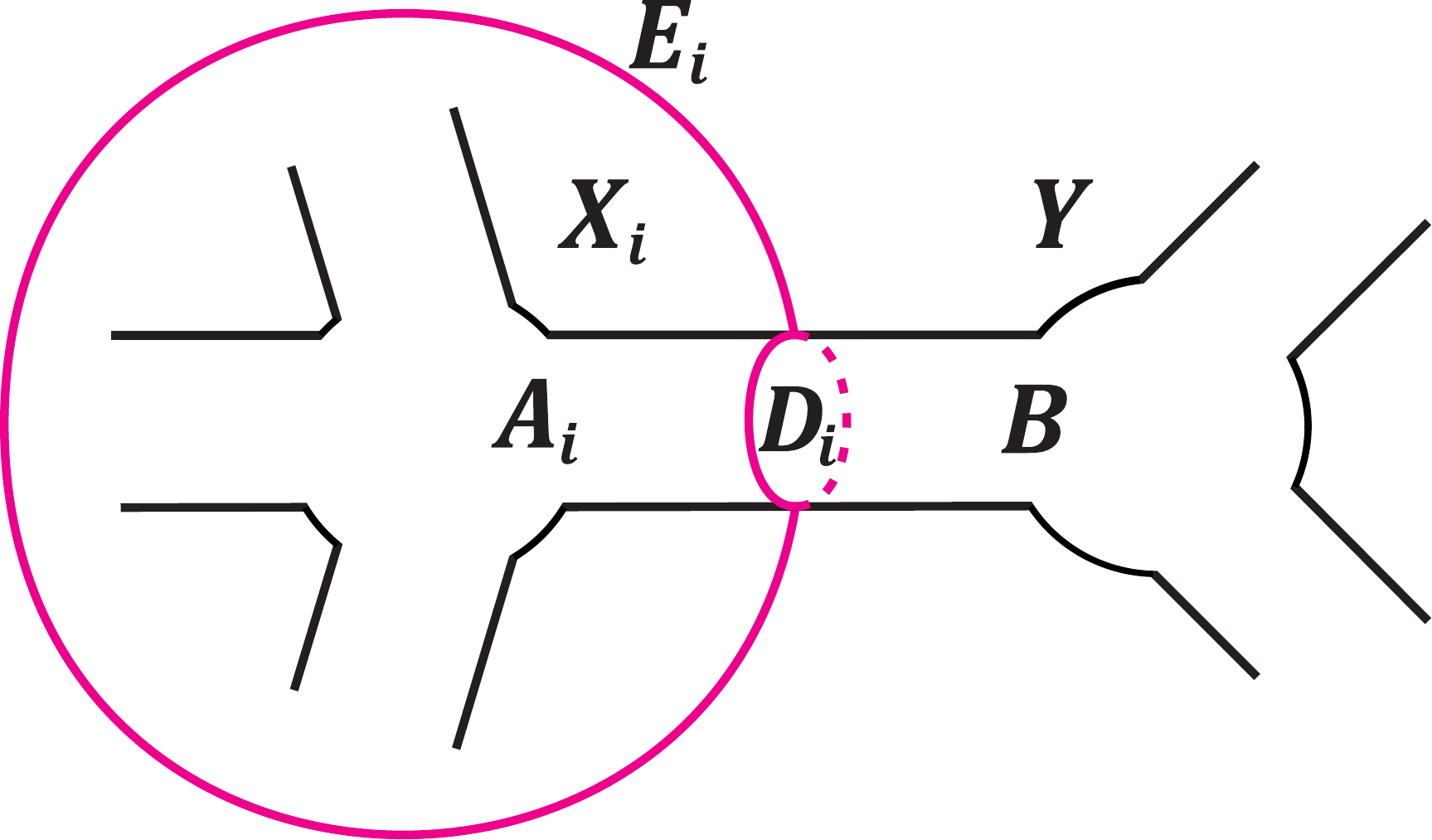}
\caption{A neighborhood $N(G)$ of $G$}\label{fig2-2}
\end{figure}

Taking a tubular neighborhood of each edge and a small ball centered at each vertex, we have a neighborhood $N(G)$ of $G$ so that it looks like Figure \ref{fig2-2}.
Let $A_i = N(G) \cap B_i$ and $B = N(G) \cap \mathrm{cl}(\mathbb{R}^3 - (\bigcup B_i))$.
Then the subset $A_i$ is a neighborhood of a union of a linearly embedded $K_{n+1}$ in $\mathrm{int}\, B_i$ and a line segment.
The subset $B$ is a neighborhood of a union of a hexagonal trefoil knot and six line segments.
It can be assumed that $D_i = N(G) \cap S_i$ is a disk, hence $E_i = \mathrm{cl} (S_i - D_i)$ is also a disk.


Let $X_i = \mathrm{cl}(B_i - A_i)$ and $Y = \mathrm{cl}(\mathbb{R}^3 - (\bigcup B_i) - B)$.
Then $X_i \cap Y =E_i$ and $X_i \cup Y= \mathrm{cl}(\mathbb{R}^3-N(G)- (\bigcup_{j \neq i} B_j))$.

\vspace{0.1cm}
\noindent {\bf Claim 1.}
{\em $\pi_1(X_i)$ is free.}

\begin{proof}
If we glue a $3$-ball $O$ to $B_i$ along their boundaries, we get the 3-sphere $\mathbb{S}^3$.
Then $X_i = \mathrm{cl}(\mathbb{S}^3-(A_i \cup O))$,
and $A_i \cup O$ can be regarded as a neighborhood of a linearly embedded $K_{n+1}$ in $\mathbb{S}^3$.
Therefore by Theorem $1$, $\pi_1(X_i)$ is free.
\end{proof}

\noindent {\bf Claim 2.}
{\em $\pi_1(Y)$ is not free.}

\begin{proof}
Note that $Y = \mathrm{cl}(\mathbb{R}^3 - (B \cup (\bigcup B_i)))$,
and $B \cup (\bigcup B_i)$ can be regarded
as a neighborhood of a hexagonal trefoil in $\mathbb{R}^3$.
Therefore $\pi_1(Y)$ is not free by the unknotting theorem in \cite{Pa}.
\end{proof}

Since $\pi_1(X_i \cap Y)$ is trivial,
$\pi_1(X_i \cup Y)$ is a free product of $\pi_1(X_i)$ and $\pi_1(Y)$ by the van Kampen theorem.
So $\pi_1(X_i \cup Y)$ is not free by Claims 1 and 2. We repeat gluing $X_{i+1}$ to $(\cdots ((Y \cup X_1)\cup X_2)\cdots)\cup X_i$. Since $Y \cup (\bigcup^{6}_{i=1} X_i) = \mathrm{cl}(\mathbb{R}^3-N(G))$, it is concluded that $\pi_{1}(\mathbb{R}^3-N(G))$ is not free.

For any nontrivial knot other than the trefoil, the construction in the above can be applied. Therefore we have an infinite family of graphs satisfying the statement of the theorem.
Note that since the number of line segments necessary to realize polygonal representation of any nontrivial knot is at least 6 \cite{Randell-2}, the number of vertices of each graph in the family is at least $6(n+1)$.

\section{Proof of Theorem \ref{thm4}}
\begin{figure}[ht!]
\centering
\includegraphics[width=7cm]{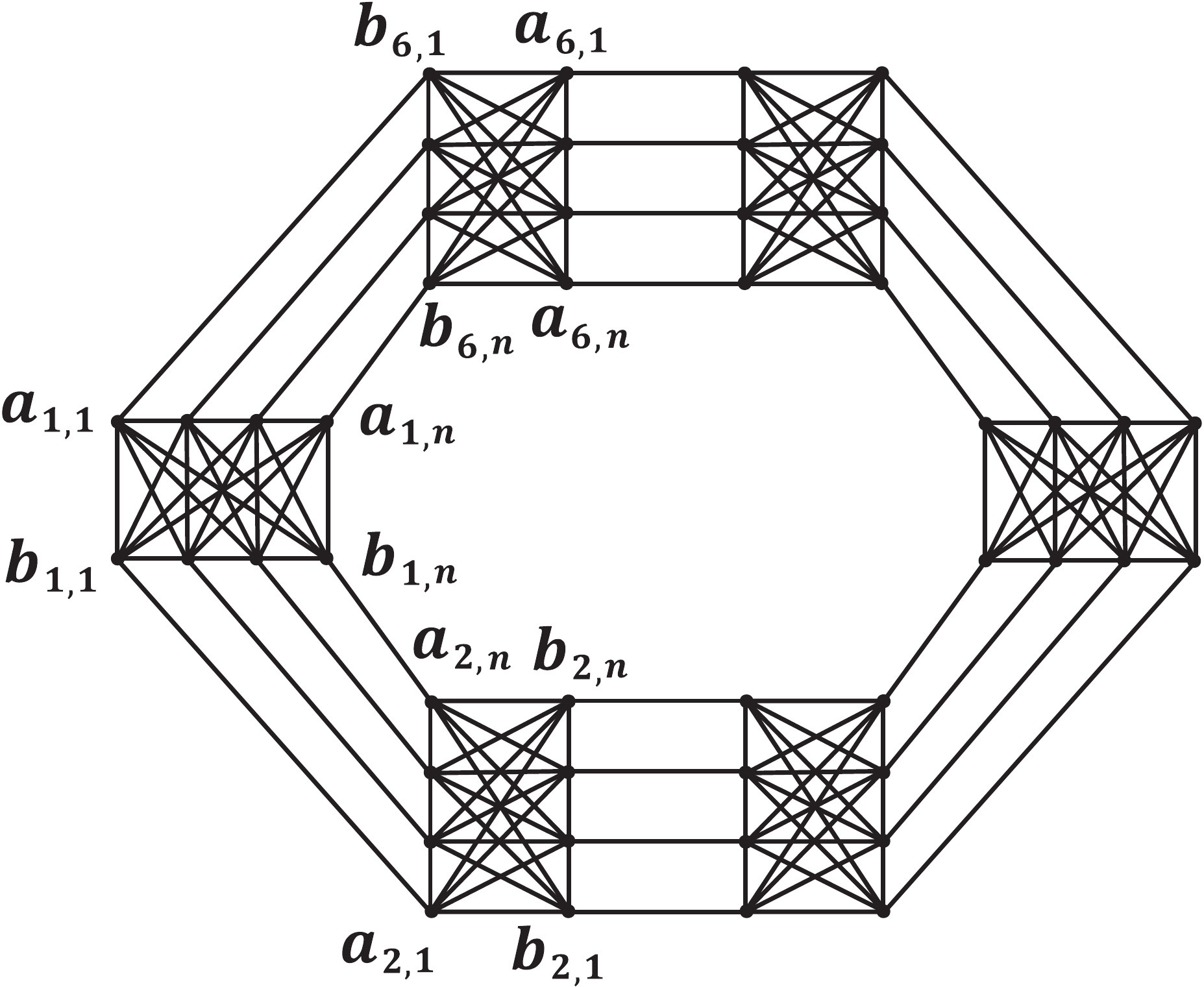}
\caption{$n$-connected graph $G$}\label{fig2-11}
\end{figure}

Let $K_{n,n}$ be the complete bipartite graph. Its vertices are denoted
by $a_i$, $b_i$ $(1 \leq i \leq n)$, and edges by $\overline{a_i b_j}$ $(1 \leq i, j \leq n)$.
And let $H$ be a graph obtained by adding edges
$\overline{a_i a_{i+1}}$ and $\overline{b_i b_{i+1}}$ $(1\leq i \leq n-1)$ to $K_{n,n}$.
Consider six copies $H_1, \ldots, H_6$ of $H$.
The vertices of each $H_i$ $(1 \leq i \leq 6)$ are denoted by
$a_{i,j}$, $b_{i,j}$ $(1 \leq j \leq n)$.
Finally we obtain a graph $G$ from $H_1, \ldots, H_6$ by adding edges
$\overline{b_{1,j} a_{2,j}}$, $\ldots$, $\overline{b_{5,j} a_{6,j}}$, $\overline{b_{6,j} a_{1,j}}$
($1 \leq j \leq n$) as illustrated in Figure \ref{fig2-11}.
Then the graph $G$ is $n$-connected.

\begin{figure}[ht!]
\centering
\includegraphics[width=7cm]{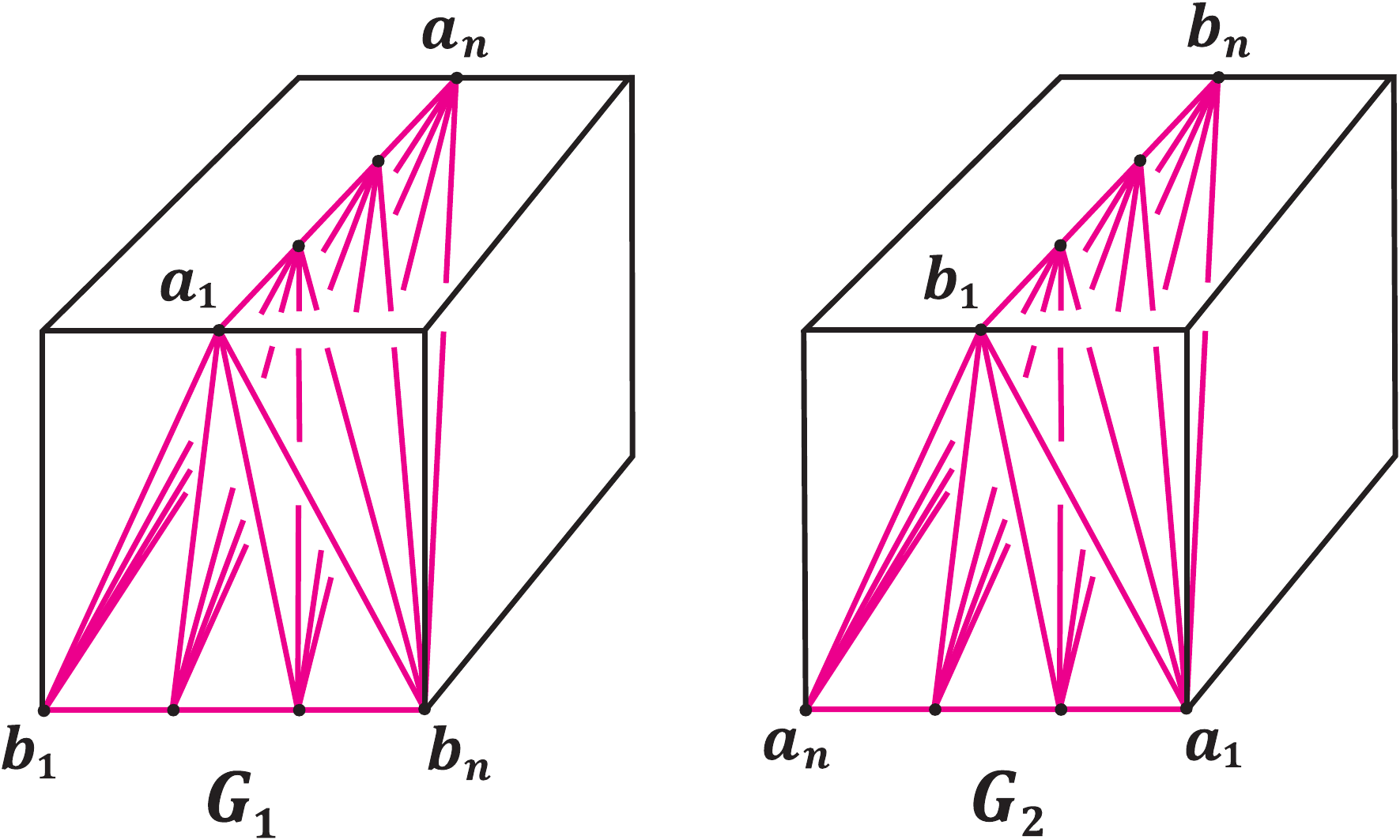}
\caption{}\label{fig2-12}
\end{figure}

Before constructing a linear embedding of $G$ into $\mathbb{R}^3$,
let us consider a linear embedding of $H$ into the cube $I^3 =[1,n] \times [1,n] \times [1,n]$.
A linear embedding of a graph is determined by the positions of its vertices.
Construct two linear embeddings $G_1$ and $G_2$ of $H$ into $I^3$ so that
$a_i = (n/2,i,n)$, $b_i = (i,1,1)$ for $G_1$, and
$a_i = (n+1-i,1,1)$, $b_i = (n/2,i,n)$ for $G_2$ as illustrated in Figure \ref{fig2-12}.

\begin{figure}[ht!]
\centering
\includegraphics[width=7cm]{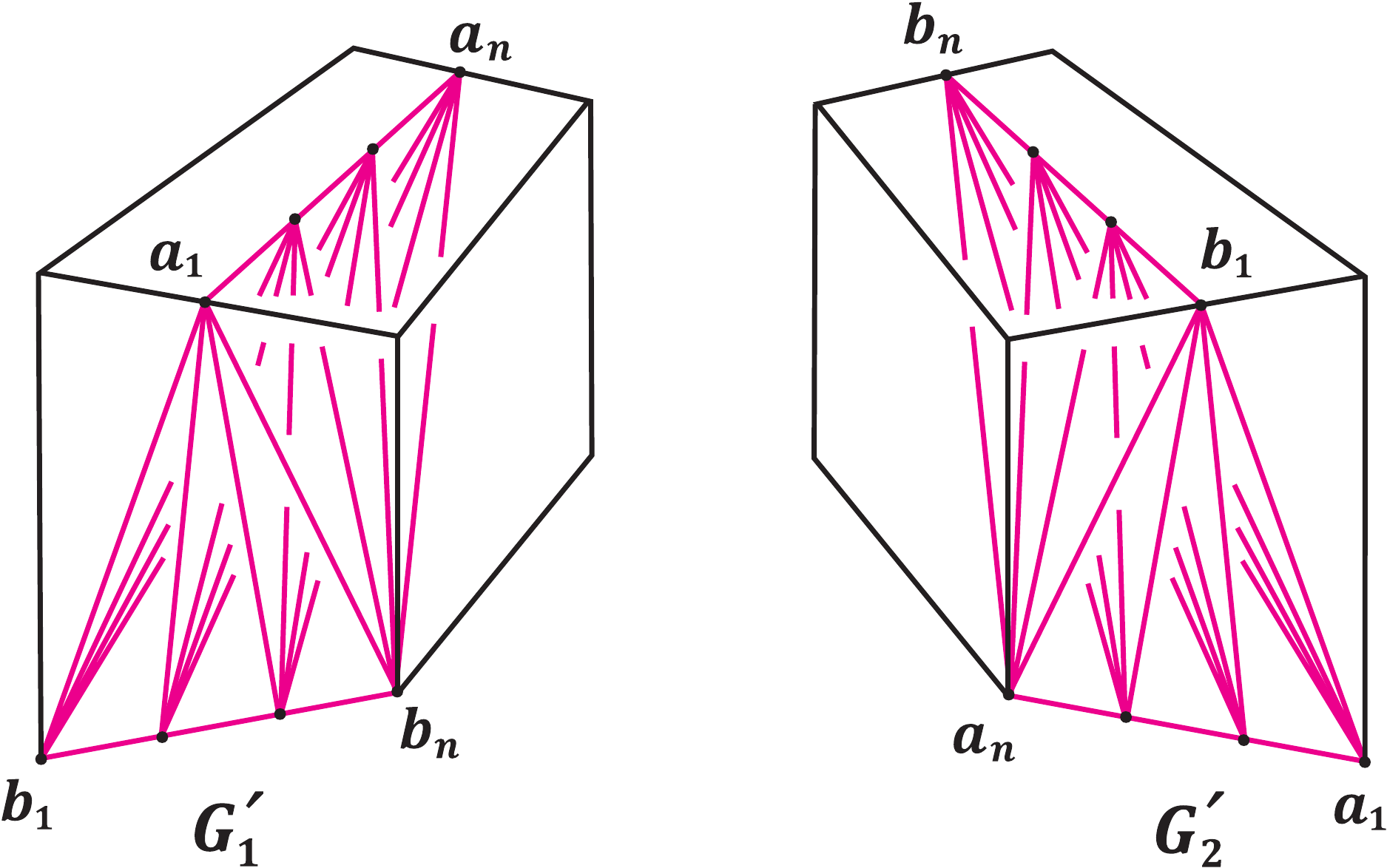}
\caption{}\label{fig2-13}
\end{figure}
\begin{figure}[ht!]
\centering
\includegraphics[width=8cm]{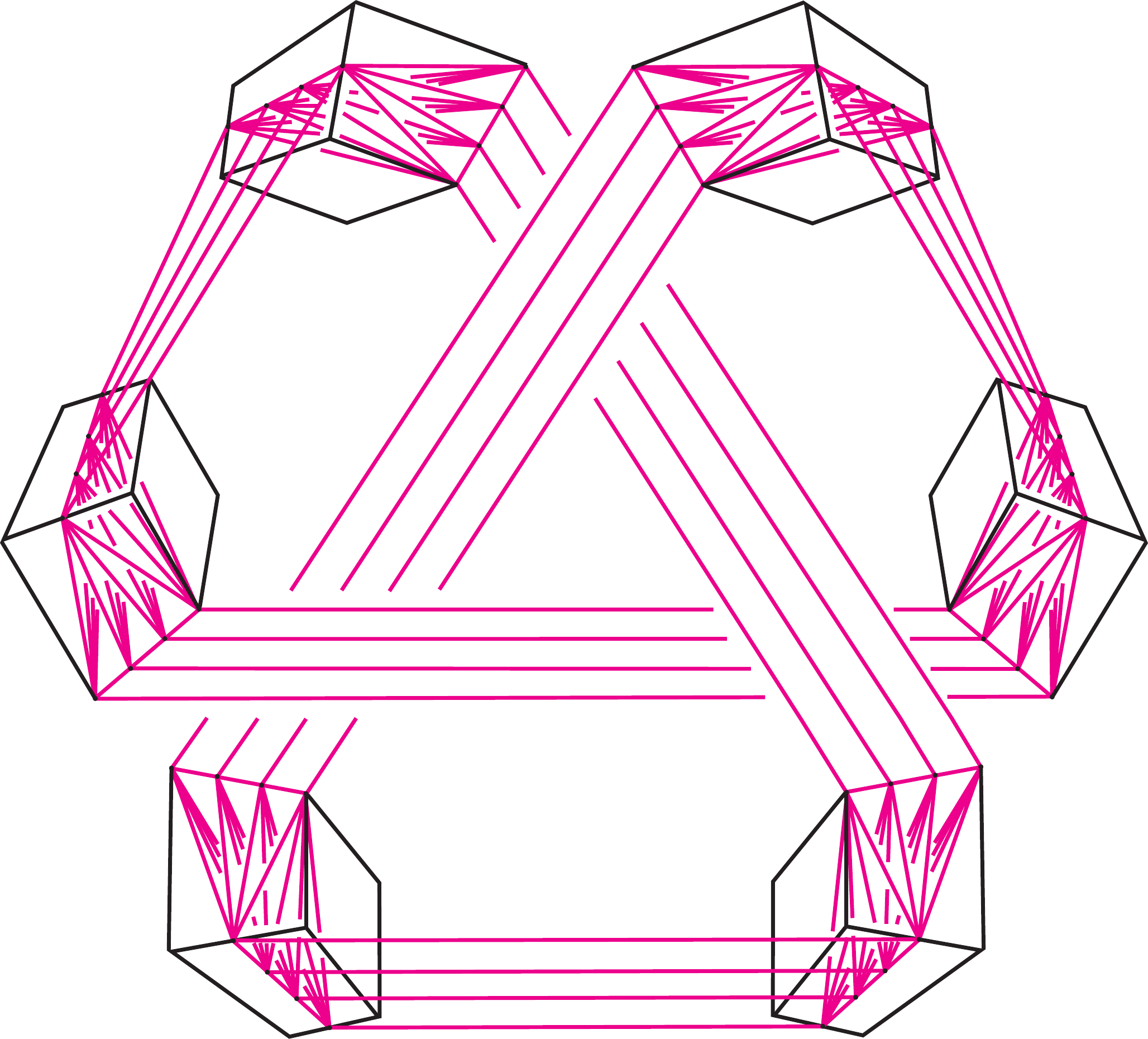}
\caption{}\label{fig2-14}
\end{figure}

Figure \ref{fig2-14} shows the final linear embedding of $G$ obtained by connecting the cubes in the above with $6n$ line segments. The embedding is realizable by modifying the cubes so that the top and front faces are leaning as much as necessary. Figure \ref{fig2-13} illustrates such modification.
Note that the embedded cycle $\langle a_{1,1} b_{1,1} a_{2,1} b_{2,1}\cdots a_{6,1} b_{6,1} \rangle$ of $G$ is a dodecagonal trefoil knot.

Consider a process as the following. Let $T_0$ be an embedded graph in $\mathbb{R}^3$. For $i \geq 0$ choose an arc $\alpha_i$ in $\mathbb{R}^3$ and a subarc $\beta_i$ of $T_i$ so that there exists a disk $D_i$
with $\partial D_i = \alpha_i \cup \beta_i$ and $D_i \cap T_i = \beta_i$. And let $T_{i+1}=T_i \cup \alpha_i$.
Then by the HNN-Extension theorem \cite{Lyndon} $\pi_1(\mathbb{R}^3-T_{i+1})$ is the free product of $\pi_1(\mathbb{R}^3-T_{i})$ and the infinite cyclic group $\mathbb{Z}$. Therefore for any $n \geq 1$ $\pi_1(\mathbb{R}^3-T_{n})$ is a free product of $\pi_1(\mathbb{R}^3-T_{0})$ and the free group of rank $n$.

We can see that the embedded graph $G$ is obtained through the process as the above starting from a dodecagonal trefoil knot. So $\pi_1(\mathbb{R}^3-G)$ is a free product of the trefoil knot group and a free group, consequently not free. By applying this construction to other nontrivial knots, we have an infinite family of graphs satisfying the statement of the theorem. Every graph in the family has at least $12n$ vertices.

\section{Proof of Theorem \ref{thm2}}
Since the minimal valency of $G$ is at least $3$, the number of vertices should be $4$, $5$ or $6$.
If it is $4$, then $G$ is the complete graph $K_4$.
By Theorem \ref{thm1} we may assume that $G$ is not a complete graph.

\begin{figure}
\centering
\includegraphics[width=11cm]{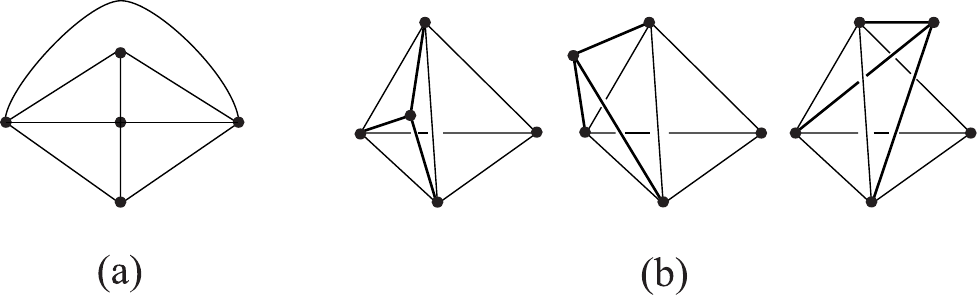}
\caption{}
\label{fig3-1-0}
\end{figure}

If $|V(G)|$ is $5$, the possible distribution of valencies of its vertices is $(4,4,4,3,3)$ or $(4,3,3,3,3)$. In the former case $G$ should be the graph in Figure \ref{fig3-1-0}-(a).
The graph contains the complete graph $K_4$ as a subgraph. If $G$ is linearly embedded into $\mathbb{R}^3$, then $K_4$ constitutes the 1-skeleton of a tetrahedron.
Consider the relative position of three edges incident to the fifth vertex $v$ with respect to the tetrahedron $T$.
Then we know that any linear embedding of $G$ corresponds to one of the three types illustrated in Figure \ref{fig3-1-0}-(b): $v \in \mbox{int}\;T$, $v \notin T$ and no edge incident to $v$ meets $\mbox{int}\; T$, or $v \notin T$ and only one edge from $v$ meets $\mbox{int}\;T$. Apply the process in the previous section which starts from $K_4$. Then it can be confirmed that the three types are free. In fact each type can be isotoped into a plane in $\mathbb{R}^3$.

\begin{figure}
\centering
\includegraphics[width=9cm]{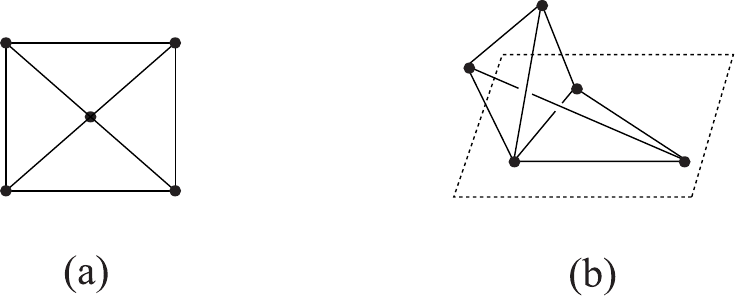}
\caption{}
\label{fig3-1}
\end{figure}

If the distribution is $(4,3,3,3,3)$, $G$ should be the graph in Figure \ref{fig3-1}-(a). Consider a convex hull determined by four vertices of $G$. Then, similarly as above, it can be shown that all possible types of linear embedding of $G$ are free. In fact, also in this case, each type can be isotoped into a plane in $\mathbb{R}^3$. Figure \ref{fig3-1}-(b) shows an example.


\vspace{2mm}
In the rest of this section we consider the case $|V(G)|=6$. Let $f$ be a linear embedding of $G$ into $\mathbb{R}^3$. For our convenience $G$ will denote both the abstract graph itself and the embedded graph without distinction. The vertices will be labelled simply by $1, 2, \ldots, 6$. And the line segment in $\mathbb{R}^3$ between two vertices $i$ and $j$ is denoted by $\overline{ij}$. Let $\Delta_{ijk}$ denote the convex hull determined by three vertices $i$, $j$ and $k$. Lastly, for an ordered sequence $ijk$ of three vertices, define
$$ H^+_{ijk} = \{p\in \mathbb{R}^3 \; | \; (\overrightarrow{ij}\times \overrightarrow{jk})\cdot \overrightarrow{jp} > 0 \} \;\; \mbox{and} \;\;
H^-_{ijk} = \{q\in \mathbb{R}^3 \; | \; (\overrightarrow{ij}\times \overrightarrow{jk})\cdot \overrightarrow{jq} < 0 \} .$$

Now three lemmas necessary for the proof of Theorem \ref{thm2} are introduced. We begin with a well-known result of Conway and Gordon.
\begin{lem} {\rm  \cite{CG}} \label{lem3-1}
Every embedding of the complete graph $K_6$ into $\mathbb{R}^3$ contains a non-splittable  $2$-component link as a pair of disjoint cycles.
\end{lem}

\begin{figure}
\centering
\includegraphics[width=7cm]{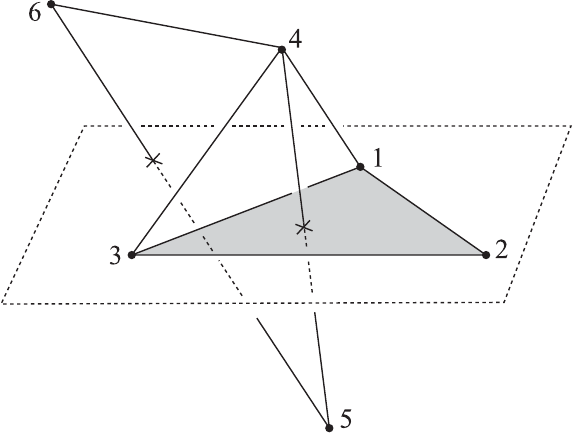}
\caption{}
\label{fig3-2}
\end{figure}

Without loss of generality we may assume that the vertices embedded by $f$ are in general position. Then $f$ can be extended to be a linear embedding of $K_6$. By Lemma \ref{lem3-1} we can label the vertices so that the following conditions are satisfied (see Figure \ref{fig3-2} for your understanding).
\begin{itemize}
\item[-] $\partial \triangle_{123} \cup \partial \triangle_{456}$
is a Hopf link.
\item[-] $\overline{45}$ penetrates $\triangle_{123}$.
\item[-] $\overline{13}$ penetrates $\triangle_{456}$.
\item[-] $4, 6 \in H_{132}^+$ and $5 \in H_{132}^-$.
\item[-] $6 \in H_{134}^+$ and $2, 5 \in H_{134}^-$.
\end{itemize}

We will say that a convex hull $\Delta_{ijk}$ is {\em trivial}, if its interior is not penetrated by any line segment between the vertices of $G$. In a previous work of the first author the linear embeddings of $K_6$ were investigated. An observation from the work is given here as a lemma.
\begin{lem} {\rm \cite[Section 4]{HJ}} \label{lem3-2}
The following convex hulls are trivial:
$$
\{ \Delta_{124}, \Delta_{125}, \Delta_{135}, \Delta_{136}, \Delta_{146}, \Delta_{156},
\Delta_{234}, \Delta_{235}, \Delta_{245}, \Delta_{246}, \Delta_{346}, \Delta_{356} \}
$$
\end{lem}

Considering the conditions of $G$, the following lemma is easily proved.
\begin{lem} \label{lem3-3}
The graph $G$ contains a Hamiltonian cycle.
\end{lem}

By Lemma \ref{lem3-3} there exists a hexagonal knot $P$ in $\mathbb{R}^3$ which is a Hamiltonian cycle of $G$. And each component of $G-P$ is a line segment connecting two vertices of $P$. Such line segments will be called {\em bridges} of $P$. Possibly $P$ is one of $\frac{5!}{2}=60$ Hamitonian cycles of $F(K_6)$, where $F$ is the linear embedding of $K_6$ extended from $f$. For each possible case,  we are going to observe the isotopy relation among the bridges of $P$.

\begin{figure}
\centering
\includegraphics[width=8cm]{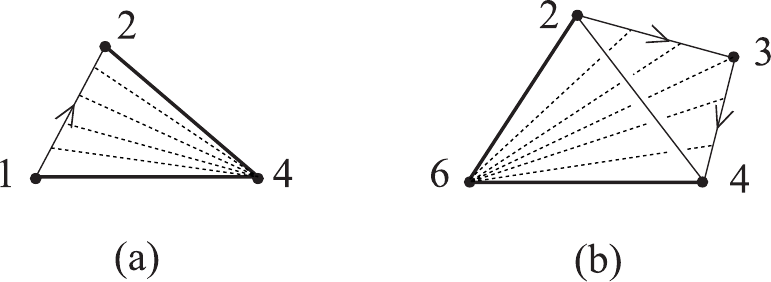}
\caption{}
\label{fig3-3}
\end{figure}

\begin{figure}
\centering
\includegraphics[width=9cm]{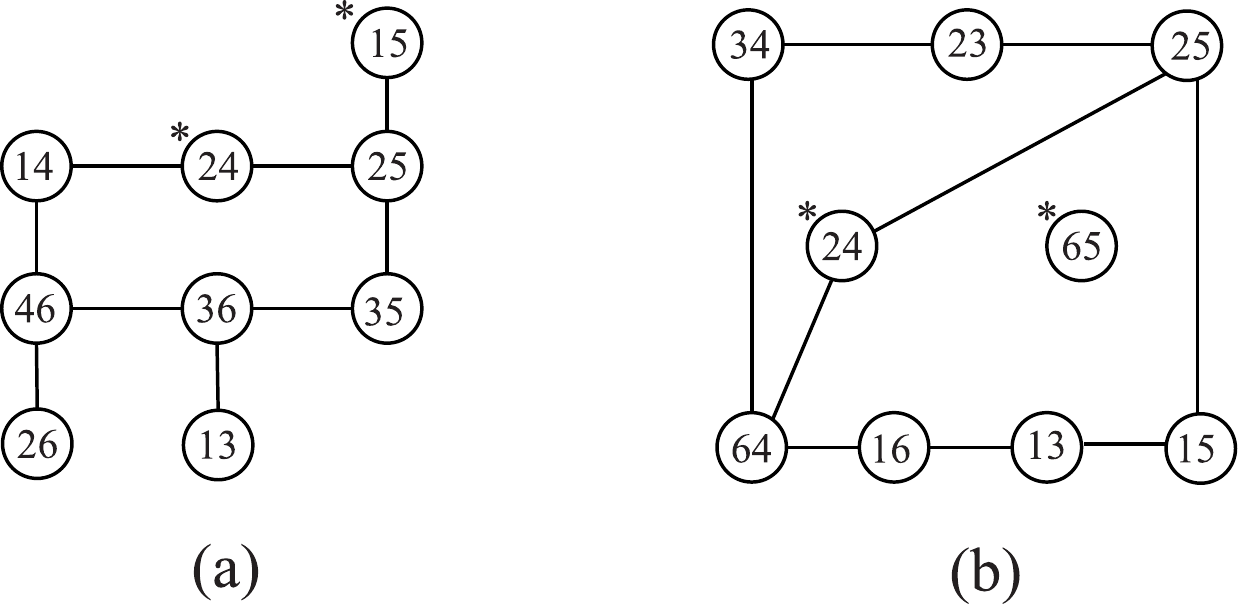}
\caption{}
\label{fig3-4}
\end{figure}

\vspace{0.5cm}
\noindent{\em Case 1: $P$ is the cycle $<123456>$.}
In this case the possible candidates for bridges of $P$ are $$\{\overline{13}, \overline{14}, \overline{15}, \overline{24}, \overline{25}, \overline{26}, \overline{35}, \overline{36}, \overline{46}\} \;.$$
Since $\Delta_{124}$ is trivial, it is a disk whose interior is disjoint from the embedded graph $K_6$. Therefore $\overline{14}$ can be isotoped to $\overline{24}$ by sliding along $\overline{12}$ (see Figure \ref{fig3-3}-(a)). Now look into $\overline{26}$ and $\overline{46}$. The trivial convex hull $\Delta_{246}$ does not contain any edge of $P$. But $\Delta_{234}$ is trivial. Therefore, slightly pulling down the disk $\Delta_{246}\cup\Delta_{234}$, we can obtain another disk such that it is bounded by $\overline{26}$, $\overline{46}$, $\overline{34}$ and $\overline{23}$, and its interior is disjoint from the embedded $K_6$. Consequently $\overline{26}$ can be isotoped to $\overline{46}$ by sliding along a subarc $\overline{23}\cup\overline{34}$ of $P$ (see Figure \ref{fig3-3}-(b)). Applying these two arguments to other candidates, we can construct a graph $I$ as seen in Figure \ref{fig3-4}-(a): The fat vertices correspond to candidates of bridges. And two fat vertices are connected by an edge, if two corresponding candidates are isotopic as described in the above. Note that the resulting graph $I$ is connected.

Now we associate the connectivity of $I$ with the freeness of $f$. Let $B_1, \ldots , B_k$ be the bridges of $P$, and $B_{k+1}, \ldots, B_{9}$ be the other line segments not belonging to $G$. Choose labels so that for $i \geq j \geq k+1$
$$d(B_i, \{B_1, \ldots, B_k\}) \geq d(B_j, \{B_1, \ldots, B_k\})\;,$$
where $d$ is the distance between vertices in the graph $I$. Let $G_0=P$ and $G_i= G_{i-1}\cup B_i$ for $i \geq 1$. And let $M_i=\mathbb{R}^3-N(G_i)$ for $i \geq 0$. Then $G_k=G$ and $G_9=K_6$.
Since the graph $I$ is connected, for any $i \geq k$, the line segment $B_{i+1}$ is parallel to $\partial M_i$. This implies that $\pi_1(M_{i+1})\cong \pi_1(M_i)*\mathbb{Z}$.
In consequence $\pi_1(M_9)$ is the free product of $\pi_1(M_k)$ and the free group of rank $9-k$.
By Theorem \ref{thm1} $\pi_1(M_9)$ is free, which implies that also $\pi_1(M_k) \cong \pi_1(\mathbb{R}^3- G)$ is free.

\vspace{0.5cm}
\noindent{\em Case 2: $P$ is the cycle $<126354>$.}
The possible candidates for bridges of $P$ are $$\{\overline{16}, \overline{13}, \overline{15}, \overline{23}, \overline{25}, \overline{24}, \overline{65}, \overline{64}, \overline{34}\} \;.$$
Figure \ref{fig3-4}-(b) shows that the graph $I$ is not connected in this case. But the line segment $\overline{65}$ bounds the trivial $\Delta_{356}$ together with $\overline{63}\cup\overline{35}$ which constitute an arc of $P$. This implies that if $\overline{65}$ is $B_i$, then it is parallel to $\partial M_{i-1}$. Also $\overline{24}$ belonging to the other component of $I$ bounds the trivial $\Delta_{124}$ together with an arc $\overline{21}\cup\overline{14}$ of $P$. Therefore the arguments in Case 1 is still valid, and we can conclude that $\pi_1(M_k)$ is free.

\vspace{0.5cm}
We had checked the remaining 58 cases (in fact, permuting the vertices $1$ and $3$, it is enough to check 28 cases). In each of the cases, either the graph $I$ is connected or every connected component of $I$ contains a fat vertex whose corresponding line segment bounds a trivial disk together with an arc of $P$. Therefore we can apply the arguments in Cases 1 and 2, to conclude that $\pi_1(M_k)$ is free.
This completes the proof.

\section*{Acknowledgements}
The first author was supported by Basic Science Research Program through the National Research Foundation of Korea(NRF) funded by the Ministry of Education, Science and Technology (2010-0009794).

The second author was supported by the National Research Foundation of Korea(NRF) grants 2011-0027989 and 2013R1A1A2059197.



\end{document}